\documentclass{amsart}
\usepackage{amssymb, amsmath, amsfonts, amscd}

\input xypic

\theoremstyle{plain}
\newtheorem{theorem}{Theorem}[section]
\newtheorem{corollary}[theorem]{Corollary}
\newtheorem{lemma}[theorem]{Lemma}
\newtheorem{proposition}[theorem]{Proposition}

\newtheorem{example}[theorem]{Example}

\theoremstyle{definition}
\newtheorem{definition}[theorem]{Definition}

\theoremstyle{remark}

\numberwithin{equation}{theorem}

\newcommand{\m}{\mathfrak{m}}
\newcommand{\F}{\mathcal{F}}

\newcommand{\I}{\mathcal{I}}

\renewcommand{\L}{\mathcal{L}}
\newcommand{\E}{\mathcal{E}}

\newcommand{\Q}{\mathcal{Q}} 
\renewcommand{\O}{\mathcal{O} }

\renewcommand{\Pr}{\mathcal{P} }

\newcommand{\SL}{\operatorname{SL}} 
\newcommand{\GL}{\operatorname{GL}} 
\newcommand{\Pic}{\operatorname{Pic} }

\newcommand{\Hom}{\operatorname{Hom} }

\newcommand{\Spec}{\operatorname{Spec} }

\renewcommand{\H}{\operatorname{H} }
\newcommand{\Proj}{\operatorname{Proj} }

\newcommand{\R}{\operatorname{R} }
\newcommand{\U}{\operatorname{U}}

\renewcommand{\lg}{\mathfrak{g}}
\newcommand{\lp}{\mathfrak{p}} 
\newcommand{\lpl}{\mathfrak{p}_L}

\renewcommand{\sl}{\mathfrak{sl}}
\renewcommand{\ln}{\mathfrak{n}}

\newcommand{\lh}{\mathfrak{h}}

\newcommand{\Pss}{P_{semi}}

\newcommand{\Z}{\mathbf{Z} }
\newcommand{\gr}{\mathbb{G} }
\newcommand{\sym}{\operatorname{Sym} }
\newcommand{\aut}{\operatorname{Aut} }

\newcommand{\ogpmod}{\underline{mod}^G(\O_{G/P})}
\newcommand{\oghmod}{\underline{mod}^G(\O_{G/H})}

\newcommand{\p}{\mathbb{P}} 

\begin{document}

\title{Jet bundles on projective space II}

\author{Helge Maakestad}

\email{\text{h\_maakestad@hotmail.com}}

\thanks{Partially supported by a scholarship from NAV, www.nav.no}


\date{January 2010}

\begin{abstract} In previous papers the structure of the jet bundle as
  $P$-module has been studied using different techniques. In this
  paper we use techniques from algebraic groups, sheaf theory, generliazed Verma
  modules, canonical filtrations of irreducible $\SL(V)$-modules and
  annihilator ideals of highest weight vectors to study the canonical
  filtration $\U_l(\lg)L^d$ of the irreducible $\SL(V)$-module 
$\H^0(X,\O_X(d))^*$ where $X=\gr(m,m+n)$.
  We study $\U_l(\lg)L^d$ using results from previous papers on the subject and recover
  a well known classification of the structure of the jet bundle $\Pr^l(\O(d))$ on
projective space $\p(V^*)$ as $P$-module. As a consequence we prove formulas
  on the splitting type of the jet bundle on projective space as
  abstract locally free sheaf. We also classify the $P$-module of the
  first order jet bundle $\Pr^1_X(\O_X(d))$ for any $d\geq 1$.
We study the incidence complex for the line bundle $\O(d)$ on the
projective line and show it is a resolution of the ideal sheaf of
$I^l(\O(d))$ - the incidence scheme of $\O(d)$.
The aim of the study is to apply it to the study of syzygies of
discriminants of linear systems on projective space and grassmannians.
\end{abstract}

\maketitle

\tableofcontents

\section{Introduction} 

In a series of papers (see \cite{flag}, \cite{maa2},\cite{maa3} and
\cite{maa4}) the structure of the jet bundle as
  $P$-module has been studied using different techniques. In this
  paper we continue this study using techniques from 
algebraic groups, sheaf theory, generalized Verma modules, canonical 
filtrations of irreducible $\SL(V)$-modules and
  annihilator ideals of highest weight vectors, and study the canonical
  filtration $\U_l(\lg)L^d$ of the $\SL(V)$-module
  $\H^0(X,\O_X(d))^*$ where $X=\gr(m,m+n)$ is the grassmannian of
  $m$-planes in an $m+n$-dimensional vector space .
  Using results obtained in \cite{flag} we classify $\U_l(\lg)L^d$ and as
  a corollary we recover a well known result on the structure of the jet bundle
  $\Pr^l(\O(d))$ on $\p(V^*)$ as $P$-module. As a consequence we get well known formulas
  on the splitting type of the jet bundle on projective space as
  abstract locally free sheaf. We also classify the $P$-module of the
  first order jet bundle $\Pr^1_X(\O_X(d))$ on any grassmannian
  $X=\gr(m,m+n)$ (see Corollary \ref{maincorr}).

In the first section of the paper we study the jet bundle
$\Pr^l_{G/H}(\E)$ of any locally free $G$-linearized sheaf $\E$ on any
quotient $G/H$. Here $G$ is an affine algebraic group of finite type
over an algebraically closed field $K$ of characteristic zero and
$H\subseteq G$ is a closed subgroup. There is an equivalence of
categories between the category of finite dimensional $H$-modules and
the category of finite rank locally free $\O_{G/H}$-modules with a $G$-linearization.
The main result of this section
is Theorem \ref{quot} where we give a classification of the $H_{l}$-modules structure of the fiber
$\Pr^l_{G/H}(\E)(x)^*$ where $H_l\subseteq H$ is a Levi
subgroup. Here $x\in G/H$ is the distinguished $K$-rational point
defined by the identity $e\in G$. We also study the structure of
$\Pr^l_X(\O_X(d))(x)^*$ as $H_l$-module where $X=\gr(m,m+n)$ is the
grassmannian of $m$-planes in an $m+n$-dimensional vector space (see
Corollary \ref{corr1} and \ref{corr2}).

In the second section we study the canonical filtration $\U_l(\lg)L^d$
for the irreducible $\SL(V)$-module $\H^0(\gr, \O_\gr(d))^*$. Here
$\gr=\gr(m,m+n)$. We prove in Theorem \ref{main} there is an
isomorphism
\[ \U_l(\lg)L^d\cong L^{d-l}\otimes \sym^l(\lg/\lpl \otimes L) \]
of $P$-modules when $\gr=\gr(1,n+1)=\p^n$ is projective $n$-space. As
a result we recover in Corollary \ref{recov} the structure of the
fiber $\Pr^l_\gr(\O_\gr(d))(x)^*$ as $P$-module. This result was proved in another paper (see
\cite{maa1}) using different techniques. We also recover in Corollary
\ref{split} a known formula on the structure of the jet bundle on
projective space as abstract locally free sheaf
(see \cite{maa2}, \cite{maa4},\cite{perk},\cite{pien},\cite{diro} and \cite{somm}).

In the third section we study the incidence complex 
\[ \wedge^\bullet \O_{\p(W^*)}(-1)_Y\otimes \Pr^l(\O(d))^*_Y \]
of the line bundle $\O(d)$ on the projective line. Using Koszul
complexes and general properties of jet bundles we prove it is a
locally free resolution of the ideal sheaf of $I^l(\O(d))$ - the
incidence scheme of $\O(d)$.

In Appendix A and B we study $\SL(V)$-modules, automorphisms of
$\SL(V)$-modules and give an elementary proof of the Cauchy formula.

Hence the paper initiates a general study of the canonical filtration $\U_l(\lg)L^d$ for 
any line bundle $\O(d)$ with $d\geq 1$ on any grassmannian $\gr(m,m+n)$ as $P$-module. In Section 3 we show 
some of the complications arising in this study by giving explicit examples.

The study of the jet bundle $\Pr^l_X(\O_X(d))$ of a line bundle
$\O_X(d)$ on the grassmannian $X=\gr(m,m+n)$ is motivated partly by
its relationship with the discriminant $D^l(\O_X(d))$ of the line
bundle $\O_X(d)$. There is by \cite{maa10} for all $1\leq l <d$ an exact sequence of locally free
$\O_X$-modules
\[ 0\rightarrow \Q \rightarrow \H^0(X,\O_X(d))\otimes \O_X
\rightarrow \Pr^l_X(\O_x(d)) \rightarrow 0 \]
giving rise to a diagram of maps of schemes
\[
\diagram \p(\Q^*) \rto^i \dto^\pi & \p(W^*)\times X \dto^p \\
         D^l(\O_X(d)) \rto^j & \p(W^*) 
\enddiagram,
\]
where $W=\H^0(X,\O_X(d))$, $\pi$ is the restriction of the projection
map and $i,j$ are closed immersions. By definition $D^l(\O_X(d)):=\pi(\p(\Q^*))$ is the schematic image
of $\p(\Q^*)$ via $\pi$.
The $K$-rational points of
$\p(\Q^*)$ are pairs of $K$-rational points $(s,x)$ with the property
that $T^l(x)(s)=0$ in $\Pr^l_X(\O_X(d))(x)$. The scheme $\p(\Q^*)$ is
the \emph{incidence scheme} of the $l$'th Taylor morphism 
\[ T^l:\H^0(X,\O_X(d))\otimes \O_X \rightarrow \Pr^l_X(\O_x(d)) .\]
The map $\pi$ is a surjective generically finite morphism between
irreducible schemes. There is by \cite{maa10} a Koszul complex of locally free sheaves on $Y=\p(W^*)\times
X$ 
\begin{align}
&\label{koszul} 0\rightarrow \O(-r)_Y\otimes \wedge^r
\Pr^l_X(\O_X(d))_Y^*\rightarrow \cdots \rightarrow \O(-1)_Y\otimes
\Pr^l_X(\O_X(d))_Y^* \rightarrow 
\end{align}
\[ \O_Y \rightarrow \O_{\p(\Q^*)}\rightarrow 0 \]
which is a resolution of the ideal sheaf of $\p(\Q^*)$ when it is
locally generated by a regular sequence. The complex \ref{koszul}
might give information on a resolution of the ideal sheaf of $D^l(\O_X(d))$.
A resolution of the ideal sheaf of $D^l(\O_X(d))$ will give information
on  its syzygies. By \cite{maa10} the first discriminant
$D^1(\O_\p(d))$ on the projective line $\p=\p^1$ is the
\emph{classical discriminant of degree d polynomials}, hence it is a
determinantal scheme. By the results of \cite{lascoux} we get an
approach to the study of the syzygies of $D^1(\O_\p(d))$. Hence we get
two approaches to the study of syzygies of discriminants of line
bundles on projective space and grassmannians: One using Taylor maps,
incidence schemes, jet bundles and generalized Verma modules. Another one using determinantal schemes.

\section{Jet bundles on quotients}

In this section we study the jet bundle of any finite rank
$G$-linearized locally free sheaf $\E$ on the grassmannian
$G/P=\gr(m,m+n)$ as $P_l$-module, where $P_l\subseteq P$ is a maximal
linearly reductive subgroup.

Let $K$ be an algebraically closed field of characteristic zero and
let $V$ be a $K$-vector space of dimension $n$. Let $H\subseteq G
\subseteq \GL(V)$ be closed subgroups.
The following holds:
There is a quotient morphism
\begin{align} \label{qu1}
&\pi:G\rightarrow G/H
\end{align}
and $G/H$ is a smooth quasi projective scheme of finite type over $K$.
Moreover
\begin{align}\label{qu2}
&H\subseteq G\text{ is parabolic if and only if }G/H\text{ is projective.}
\end{align}
For a proof see \cite{jantzen}. Let $X=G/H$ and let $\oghmod$ be the
category of locally free $\O_{G/H}$-modules with a
$G$-linearization. Let $\underline{mod}(H)$ be the category of finite
dimensional $H$-modules.
It follows from \cite{jantzen} there is an exact equivalence of categories
\[ \underline{mod}(H)\cong  \oghmod  .\]
Let $\E \in  \oghmod$  be a locally free $\O_{G/H}$-module. 

Let $Y=G/H\times G/H$ and $p,q:Y\rightarrow G/H$ be the canonical projection maps. 
The scheme $G/H$ is smooth and separated over $\Spec(K)$ hence the
diagonal morphism
\[ \Delta:G/H\rightarrow Y \]
is a closed immersion of schemes.
Let $\I\subseteq \O_Y$ be the ideal of the diagonal and let
$\O_{\Delta^l}=\O_Y/\I^{l+1}$ be the structure sheaf of the \emph{$n$'th
infinitesimal neigborhood of the diagonal}.
\begin{definition} Let $\E$ be a locally free finite rank
  $\O_{G/H}$-module.
Let
\[ \Pr^l_{G/H}(\E)=p_*(\O_{\Delta^l}\otimes q^*\E) \]
be the $l$'th jet bundle of $\E$.
\end{definition}

\begin{proposition} There is for all $l\geq 1$ an exact sequence of
  locally free $\O_{G/H}$-modules
\begin{align} \label{fund} 0\rightarrow \sym^l(\Omega^1_{G/H})\otimes \E
  \rightarrow 
\Pr^l_{G/H}(\E)\rightarrow^\phi \Pr^{l-1}_{G/H}(\E)\rightarrow 0 
\end{align}
with $G$-linearization.
\end{proposition}
\begin{proof} By \cite{maa4} sequence \ref{fund} is an exact sequence
  of locally free $\O_{G/H}$-modules. 
The scheme $Y$ is equipped with the diagonal $G$-action. It
follows $p_*$ and $q^*$ preserve $G$-linearizations. We get a diagram
of exact sequences of $\O_Y$-modules with a $G$-linearization
\[
\diagram 0\rto & \I^{l+1}\otimes q^*\E \rto \dto & \O_Y \otimes q^*\E
\rto \dto & \O_{\Delta^l}\otimes q^*\E \rto \dto & 0 \\
0 \rto & \I^{l}\otimes q^*\E \rto & \O_Y \otimes q^*\E \rto &
\O_{\Delta^{l-1}}\otimes q^*\E \rto & 0
\enddiagram.
\]
Since $p_*$ preserves $G$-linearization we get a morphism
\[ \phi: \Pr^l_{G/H}(\E)\rightarrow \Pr^{l-1}_{G/H}(\E) \]
preserving the $G$-linearization, and the Proposition is proved.
\end{proof}

Let $\lg=Lie(G)$ and $\lh=Lie(H)$.  Let $H_l\subseteq H$ be a Levi
subgroup of $H$. It follows $H_l$ is a maximal linearly reductive
subgroup of $H$. The group $H_l$ is not unique but all such groups are
conjugate under automorphisms of $H$. Let $x\in G/H$ be the
$K$-rational point defined by the identity $e\in G$.

\begin{theorem} \label{quot} There is for all $l\geq 1$ an isomorphism
\begin{align}\label{iso}
\Pr^l_X(\E)(x)^*\cong \E(x)^*\otimes (\oplus_{i=}^l\sym^i(\lg/\lh))
\end{align}
of $L$-modules.
\end{theorem}
\begin{proof} Dualize the sequence \ref{fund} and take the fiber at
  $x$ to get the exact sequence
\[ 0\rightarrow \Pr^{l-1}_X(\E)(x)^*\rightarrow
\Pr^l_X(\E)(x)^*\rightarrow \E(x)^*\otimes \sym^l(\lg/\lh) \rightarrow
0 \]
of $H$-modules (and $H_l$-modules). This sequence splits since $H_l$ is
linearly reductive and the Theorem follows by induction on $l$.
\end{proof}

Hence the study $\Pr^l_X(\E)(x)^*$ as $H_l$-module is reduced to the
study  of $\E(x)^*$ and $\sym^l(\lg/\lh)$.

Let $W\subseteq V$ be $K$-vector spaces of dimension $m$ and $m+n$
and let $G=\SL(V)$ and $P\subseteq G$ the subgroup fixing $W$. It
follows $G/P=\gr(m,m+n)$ is the grassmannian of $m$-planes in $V$. Let
$\lg=Lie(G)$ and $\lp=Lie(P)$. Fix a basis $e_1,..,e_m$ for $W$ and
$e_1,..,e_m,e_{m+1},..,e_{m+n}$ for $V$. It follows the $K$-rational
points of $P$ are matrices $M$ on the form

\[
M=\begin{pmatrix} A & X \\
                0 & B 
\end{pmatrix}
\]
where $det(A)det(B)=1$, $A$ an $m\times m$-matrix and $B$ an
$n\times n$-matrix. Let $P_l\subseteq P$ be the subgroup defined as
follows: The $K$-rational points of $P_l$ are matrices $M$ on the form
\[M=
\begin{pmatrix} A & 0 \\
                0 & B
\end{pmatrix}
\]
where $det(A)det(B)=1$ and similarly $A$ an $m\times m$-matrix and $B$
an $n\times n$-matrix. It follows $P_l$ is a Levi subgroup of $P$,
hence it is a maximal linearly reductive subgroup.

\begin{proposition} \label{canonical} There is a canonical isomorphism
\[ \lg/\lp\cong \Hom(W,V/W) \]
of $P$-modules.
\end{proposition}
\begin{proof} By definition $\lg=\sl(V)$, hence $\phi\in \lg$ is a map
\[ \phi: V\rightarrow V \]
with $tr(\phi)=0$. Let $i:W\rightarrow V$ be the inclusion map and
$p:V\rightarrow V/W$ the projection map. Define the following map:
\[ \j':\lg\rightarrow \Hom(W,V/W) \]
by
\[ j'(\phi)=p\circ \phi \circ i .\]
It follows $j(\lp)=0$ hence we get a well defined map
\[ j:\lg/\lp\rightarrow \Hom(W,V/W) \]
defined by 
\[ j(\overline{\phi})=p\circ\phi \circ i.\]
One checks $\lg/\lp$ and $\Hom(W,V/W)$ are $P$-modules and $j$ a
morphism of $P$-modules. It is an isomorphism and the Proposition follows.
\end{proof}

\begin{corollary} \label{corr1} On $X=\gr(m,m+n)$ there is an isomorphism
\[ \Pr^l_X(\E)(x)^*\cong \E(x)^*\otimes (\oplus_{i=0}^l
\sym^i(\Hom(W,V/W)) \]
of $P_l$-modules.
\end{corollary}
\begin{proof} The proof follows from Theorem \ref{quot} and
  Proposition \ref{canonical}.
\end{proof}

There is an isomorphism of $P$-modules 
\[ \Hom(W,V/W)\cong W^*\otimes V/W \]
hence the decomposition into irreducible components of the module
$\sym^i(W^*\otimes V/W)$ as $P_l$-module may be done using the \emph{Cauchy formula}
(see Appendix B).

Let $\lambda - |i|$ denote $\lambda$ is a  partition of the integer $i$
If $\lambda=\{ \lambda_1,..,\lambda_d\}$ is a partition of an integer
$l$,
let $\mu(\lambda)$ denote the following partition:
\[ \mu(\lambda)_i=l-\lambda_{d+1-i}.\]
Let for any partition $\lambda$ of an integer $l$ and any vector space
$W$, $\mathbb{S}_\lambda(W)$ denote the \emph{Schur-Weyl module} of
$\lambda$.

\begin{corollary} There is an isomorphism
\[ \Pr^l_X(\E)(x)^*\cong \E(x)^*\otimes
(\oplus_{i=0}^l(\bigoplus_{\lambda-|i|}\mathbb{S}_\lambda(W^*)\otimes
\mathbb{S}_{\mu(\lambda)}(V/W))) \]
of $\SL(W)\times \SL(V/W)$-modules.
\end{corollary}
\begin{proof} By Corollary \ref{corr1} there is an isomorphism
\[ \Pr^l_X(\E)(x)^*\cong \E(x)^*\otimes (\oplus_{i=0}^l
\sym^i(\Hom(W,V/W)) \]
of $P_l$-modules and $\SL(W)\times \SL(V/W)$-modules, since
$\SL(W)\times \SL(V/W)\subseteq P_l$ is a closed subgroup.
Since
\[ \sym^i(\Hom(W,V/W))\cong \sym^i(W^*\otimes V/W) \]
the result follows from the Cauchy formula (see Appendix B or \cite{fulton}).
\end{proof}

\begin{example} Calculation of the cohomology group $\H^i(X,\wedge^j\Pr^l_X(\O_X(d))^*)$.
\end{example}
In the following we use the notation introduced in \cite{jantzen}.
Let $P_{semi}=\SL(m)\times \SL(n)\subseteq P$ be the semi
simplification of $P$. We get a vector bundle
\[ \pi:G/P_{semi}\rightarrow G/P=\gr(m,m+n) .\]
Let $X=G/P$ and $Y=G/P_{semi}$
Given any finite dimensional $P$-module $W$, let $\L_{X}(W)$ denote
its corresponding $\O_{X}$-module. Let $W_{semi}$ denote the
restriction of $W$ to $P_{semi}$. By the results of \cite{jantzen}
it follows there is an isomorphism
\[ \pi^*\L_{X}(W)\cong \L_{Y}(W_{semi}) \]
of locally free sheaves. This will help calculating the higher
cohomology group 
\[ \H^i(X, \L_{X}(W)) \]
since $P_{semi}$ is semi simple and $\pi$ is a locally trivial
fibration. If $W$ is the $P$-module corresponding to the dual of the $j$'th
exterior power of the jet bundle $\wedge^j \Pr^l_{X}(\O_{X}(d))^*$
we can use this construction to calculate the cohomology group
\[ \H^i(X, \wedge^j\Pr^l_{X}(\O_{X}(d))^*) .\]
Such a calculation will by the results of \cite{maa10}, Example 5.12 give
information on resolutions of the ideal sheaf of $D^l(\O_{X}(d))$
since the push down of the Koszul complex \ref{koszul} is the locally trivial sheaf
\[ \O(-j)\otimes \H^i(X, \wedge^j\Pr^l_{X}(\O_{X}(d))^*) .\]
To describe the locally trivial sheaf $\O(-j)\otimes \H^i(X, \wedge^j\Pr^l_{X}(\O_{X}(d))^*)$
for all $i,j$ we need to calculate the dimension $h^i(X,
\wedge^j\Pr^l_{X}(\O_{X}(d))^*)$
and this calculation may be done using the approach indicated above.

Let $m=2, n=4$ and $X=\gr(2,4)$. 

\begin{corollary} \label{corr2} There is an isomorphism
\[ \Pr^l_X(\E)(x)^*\cong \E(x)^*\otimes (\oplus_{i=0}^l\oplus_{j=0}^n
\sym^{2j+m}(W^*)\otimes \sym^{2j+m}(V/W)) \]
of $\SL(2)\times \SL(2)$-modules. Here $(n,m)=(\frac{i}{2},0)$ if $i=2n$ and
$(\frac{i-1}{2},1)$ if $i=2n+1$.
\end{corollary}
\begin{proof} This follows from Corollary \ref{corr1} and Proposition \ref{prop1}.
\end{proof}

\section{On canonical filtrations and jet bundles on projective space}

In this section we study the canonical filtration for the dual of the
$\SL(V)$-module of global sections of an invertible sheaf on the
grassmannian. We classify the canonical filtration on projective space
and as a result recover known formulas on the splitting type of the
jet bundle as abstract locally free sheaf.

Let $W\subseteq V$ be vector spaces over $K$ of dimension $m$ and
$m+n$.
Let $W$ have basis $e_1,..,e_m$ and $V$ have basis $e_1,..,e_{m+n}$.
Let $V^*$ have basis $x_1,..,x_{m+n}$. Let $G=\SL(V)$ and $P\subseteq
G$ the parabolic subgroup of elements fixing $W$. It follows there is
a quotient morphism
\[ \pi:G\rightarrow G/P \]
and $G/P\cong \gr(m,m+n)$ is the grassmannian of $m$-planes in $V$.
Let $\p=\gr(1,n+1)=\p(V^*)$.
Let $L^d=\sym^d(\wedge^mW)$. There is an inclusion of $P$-modules
$L^d\subseteq \sym^d(\wedge^m V)$. Since $K$ has characteristic zero
there is an inclusion of $G$-modules
\[ \H^0(G/P,\O_{G/P}(d))^*\subseteq \sym^d(\wedge^mV^*)^*\cong
\sym^d(\wedge^mV) .\]
Let $\lg=Lie(G)$ and $\lp=Lie(P)$. Let $\U(\lg)$ be the universal
enveloping algebra og $\lg$ and let $\U_l(\lg)$ be the $l$'th term ot
its canonical filtration.

By the Corollary 3.11 in \cite{maa11} there is for all $1\leq l \leq d$ an exact sequence
of $P$-modules
\[ 0\rightarrow \Pr^l_\gr(\O_\gr(d))(x)^*\rightarrow
\H^0(\gr,\O_\gr(d))^*\rightarrow \H^0(\gr,\m^{l+1}\O_\gr(d))^*\rightarrow
0  .\] 
Since the grassmannian is projectively normal in the Plucker embedding
we get an inclusion 
\[ \H^0(\gr,\O_\gr(d))^*\subseteq \sym^d(\wedge^m V)\]
of $P$-modules. The highest weight vector for $\H^0(\gr,\O_\gr(d))^*$ is the line
$L^d=\sym^d(\wedge^m W)$. 
Let $ann(L^d)\subseteq \U(\lg)$ be the left annihilator ideal of
$L^d$. It is the ideal generated by elements $x\in \U(\lg)$ with the
property $x(L^d)=0$. Let $ann_l(L^d)$ be its canonical filtration. 
We get an exact sequence of $G$-modules
\[ 0\rightarrow ann(L^d)\otimes L^d \rightarrow \U(\lg)\otimes L^d
\rightarrow \H^0(X,\O_X(d))^* \rightarrow 0\]
and an exact sequence of $P$-modules
\[ 0\rightarrow ann_l(L^d)\otimes L^d\rightarrow \U_l(\lg)\otimes L^d
\rightarrow \U_l(\lg)L^d \rightarrow 0\]
for all $l\geq 1$.
The $G$-module $\U(\lg)\otimes L^d$ is the \emph{generalized Verma
  module} corresponding to the $P$-module defined by $L^d\in \sym^d(\wedge^m V)$.
There is an inclusion of $P$-modules
\[ \U_l(\lg)L^d\subseteq \H^0(\gr,\O_\gr(d))^*. \]
\begin{definition} Let $\{\U_l(\lg)L^d\}_{l\geq 1}$ be the
  \emph{canonical filtration} for $\H^0(\gr ,\O_\gr(d))^*$.
\end{definition}

\begin{lemma} Assume $y\in \lg$ and $x_1\cdots x_i\in \U_i(\lg)$ with
  $x_i\in \lg$. The following holds:
\[ y(x_1\cdots x_i)=(x_1\cdots x_i)y+\omega \]
where $\omega\in \U_{i-1}(\lg)$.
\end{lemma}
\begin{proof} The proof is by induction.
\end{proof}

The Lie algebra $\lp$ is the sub Lie algebra of $\lg=\sl(V)$ given by
matrices $M$ of the following type:
\[
M=\begin{pmatrix} A & X \\
                  0 & B
\end{pmatrix}
\]
where $A$ is an $m\times m$-matrix, $B$ and $n\times n$-matrix and $tr(A)+tr(B)=0$.
Let $\lpl$ be the sub Lie algebra of $\lp$ consisting of matrices
$M\in \lp$
of the following type:
\[
M=\begin{pmatrix} A & X \\
                  0 & B
\end{pmatrix}
\]
where $tr(A)=tr(B)=0$.

\begin{proposition}\label{mainprop}
\begin{align} \label{no1} \text{The sub Lie algebra $\lpl \subseteq \lp$ is a
  sub $P$-module of $\lp$.}
\end{align}

There is an exact sequence of $P$-modules
\begin{align} \label{no2}
0\rightarrow \lp/\lpl \rightarrow \lg/\lpl \rightarrow \lg/\lp
\rightarrow 0
\end{align}
and $\lp/\lpl$ is the trivial $P$-module.

The following holds:
\begin{align}\label{no3}
 dim_K( L^{d-k}\otimes \sym^k(\lg/\lpl \otimes L) )=\binom{mn+k}{mn}.
\end{align}

There is a filtration of $P$-modules
\begin{align}\label{no4}
0=G_{l+1}\subseteq G_l\subseteq \cdots \subseteq G_0=L^{d-l}\otimes
\sym^l(\lg/\lpl \otimes L) 
\end{align}
with quotients
\[ G_i/G_{i+1}\cong L^{d-(l-i)}\otimes \sym^{l-i}((\lg/\lp \otimes L) \]
for $1\leq i \leq k$.

Assume $dim_K(W)=1$ and let $W=L$.\label{isoV} There is an exact
sequence of $P$-modules
\begin{align} \label{no5} 0\rightarrow \lpl \otimes L \rightarrow \lg \otimes L \rightarrow V
\rightarrow 0
\end{align}
giving an isomorphism of $P$-modules $\lg/\lpl\otimes L\cong V$.
\end{proposition}
\begin{proof} We prove \ref{no1}: In the following $A,a$ are square
  matrices of size $m$
  and $b,B$ square matrices of size $n$.
The $K$-rational points of the group $P$ are matrices
  $g$ on the form
\[
g=
\begin{pmatrix} A & X \\
                0 & B

\end{pmatrix}
\]
where $det(A)det(B)=1$. Assume $x\in \lp$ is the following element:
\[
x=
\begin{pmatrix} a & x \\
                0 & b
\end{pmatrix}
\]
with $tr(a)+tr(b)=0$.
It follows $g(x)=gxg^{-1}$ has $tr(gxg^{-1})=tr(gg^{-1}x)=tr(x)=0$
hence $gxg^{-1}\in \lp$ and $\lp$ is a $P$-module. Assume $x\in \lpl$
ie $tr(a)=tr(b)=0$. It follows

\[
gxg^{-1}=
\begin{pmatrix} aAa^{-1} & * \\
                 0 & bBb^{-1}
\end{pmatrix}
\]
and $tr(aAa^{-1})=tr(aa^{-1}A)=tr(A)=0$ hence $g(x)\in \lpl$ and
\ref{no1} is proved.

We prove \ref{no2}: By \ref{no1} it follows $\lpl \subseteq \lp$
is a sub
$P$-module. One checks $\lp/\lpl$ is a trivial $P$-module. We clearly
get an exact sequence of $P$-modules and \ref{no2} is proved.

We prove \ref{no3}: Since
\[ dim_K(\lg)=(m+n)^2-1=n^2+2mn+m^2-1 \]
and
\[ dim_K(\lpl)=m^2+mn+n^2-2 \]
it follows $dim_K(\lg/\lpl)=mn+1$.
It follows
\[ dim_K(L^{d-l}\otimes \sym^l(\lg/\lpl \otimes
L)=\binom{mn+1+l-1}{mn+1-1}=\binom{mn+l}{mn} .\]

We prove \ref{no4}: Since $\lp/\lpl$ is a trivial $P$-module there are
  isomorphisms of $P$-modules
\[ L^{d-(k-i)}\otimes \sym^{k-i}(\lg/\lpl \otimes L)\cong
L^{d-k}\otimes L^i\otimes \sym^{k-i}(\lg/\lpl \otimes L)\cong \]
\[ L^{d-k}\otimes \sym^i(\lp/\lpl \otimes L)\otimes
\sym^{k-i}(\lg/\lpl \otimes L)\]
for all $1\leq i \leq k$.
We get an injection
\[j:L^{d-k}\otimes \sym^i(\lp/\lpl \otimes L)\otimes
\sym^{k-i}(\lg/\lpl \otimes L) \rightarrow L^{d-k}\otimes
\sym^k(\lg/\lpl \otimes L) \]
defined by
\[ j(L^{d-k}\otimes \overline{y_1}\otimes L\cdots
\overline{y_i}\otimes L\otimes \overline{x_1}\otimes L\cdots
\overline{x_{k-i}}\otimes L)=L^{d-k}\otimes \overline{y_1}\otimes
L\cdots \overline{y_i}\otimes L \overline{x_1}\otimes L\cdots
\overline{x_{k-i}}.\]
The injection $j$ gives rise to an injection
\[ L^{d-(k-i)}\otimes \sym^{k-i}(\lg/\lpl \otimes
L)\cong L^{d-k}\otimes \sym^i(\lp/\lpl \otimes L)\otimes
\sym^{k-i}(\lg/\lpl \otimes L)\rightarrow^j \]
\[ L^{d-k}\otimes
\sym^k(\lg/\lpl \otimes L) \]
of $P$-modules for all $1\leq i \leq k$. The exact sequence
\[ 0\rightarrow \lp/\lpl \rightarrow \lg/\lpl \rightarrow \lg/\lp
\rightarrow 0 \]
gives rise to a filtration of $P$-modules
\[0=F_{l+1}\subseteq F_l\subseteq \cdots \subseteq F_0=\sym^l(\lg/\lpl
\otimes L) \]
with quotients
\[F_i/F_{i+1}\cong L^i\otimes \sym^{l-i}(\lg/\lp \otimes L).\]
Put $G_i=L^{d-l}\otimes F_i.$
It follows
\[ G_i=L^{d-(l-i)}\otimes \sym^{l-i}(\lg/\lpl \otimes L) .\]
There is an isomorphism
\[ G_i/G_{i+1}\cong L^{d-(l-i)}\otimes \sym^{l-i}(\lg/\lp \otimes L) \]
and claim \ref{no4} is proved.

We prove \ref{no5}: 
Let $V=K\{e_0,..,e_n\}$ and $L=W=e_0$. It follows
  $P\subseteq G=\SL(V)$ is the group whose $K$-rational points are the
  following:
\[
g=\begin{pmatrix} a & * \\
                  0 & B
\end{pmatrix}
\]
with $a=\frac{1}{det(B)}$. Also $B$ is an $n\times n$-matrix with
coefficients in $K$. By definition the maps in the sequence are maps
of $P$-modules. It follows $\lp=Lie(P)$ is the Lie algebra whose
elements $x$ are matrices on the following form:
\[
x=\begin{pmatrix} -tr(B) & * \\
                   0     & B
\end{pmatrix}
\]
where $B$ is any $n\times n$-matrix with coefficients in $K$. The sub
Lie algebra $\lpl \subseteq \lp$ is the Lie algebra of matrixes $x\in
\lp$ on the following form:
\[
x=\begin{pmatrix} 0 & * \\
                  0 & B
\end{pmatrix}
\]
where $B$ is any $n\times n$-matrix with $tr(B)=0$.
Let $x_i\in \lg$ be the following element: Let the first column vector
of $x_i$ be the vector $e_i$ and let the rest of the entries be such
that $tr(x_i)=0$.
It follows $x_i\otimes e_0\in \lg\otimes L$ and $x_i(e_0)=e_i$ hence
the vertical map is surjective. One easily checks the sequence is
exact and \ref{no5} is proved.
\end{proof}

We get two $P$-modules: $\lpl \subseteq \lp$ and $L^i=\sym^i(\wedge^m W)\subseteq
\sym^i(\wedge^m V)$. We get for all $1\leq k \leq d$ a $P$-module
\[ L^{d-k}\otimes \sym^k(\lg/\lpl \otimes L) .\]
There is an injection of $P$-modules
\[ i: L^{d-k}\otimes \sym^k(\lg/\lpl\otimes L)\rightarrow
\sym^d(\wedge^m V) \]
defined by
\[ i(L^{d-k}\otimes \overline{x_1}\otimes L\cdots
\overline{x_k}\otimes
L)=L^{d-k}x_1(L)\cdots x_k(L).\]

There are natural embeddings of $P$-modules
\[ \U_k(\lg)L^d \subseteq \sym^d(\wedge^m V) \]
and
\[L^{d-(k-1)}\otimes \sym^{k-1}(\lg/\lpl\otimes L)\subseteq 
 L^{d-k}\otimes \sym^k(\lg/\lpl\otimes L) \subseteq \sym^d(\wedge^m V).\]

Assume in the following $m=1$ and $L=W$. It follows $\gr=\p(V^*)=\p$ is
projective $n$-space.
\begin{proposition} \label{envelop} Let $x_1\cdots x_k(L^d)\in \U_k(\lg)L^d$. The
  following formula holds:
\[x_1\cdots x_k(L^d)=\alpha L^{d-k}x_1(L)\cdots x_k(L)+\omega \]
where $\omega\in L^{d-(k-1)}\otimes \sym^{k-1}(\lg/\lpl \otimes L)$.
\end{proposition}
\begin{proof} we prove the result by induction on $k$. Assume $k=1$ 
and let $x(L^d)\in \U_1(\lg)L^d$. It follows $x(L^d)=dL^{d-1}x(L)\in
L^{d-1}\otimes \sym^1(\lg/\lpl \otimes L)$ and the claim holds for
$k=1$. Assume the result is true for $k$. Hence
\[x_1\cdots x_k(L^d)=\alpha L^{d-k}x_1(L)\cdots x_k(L)+\omega \]
with $\omega \in L^{d-(k-1)}\otimes \sym^{k-1}(\lg/\lpl \otimes L)$. 
Assume 
\[ \omega=\sum_i \alpha_i L^{d-(k-1)}x^i_1(L)\cdots x^i_{k-1}(L) .\]
We get
\[x_0x_1\cdots x_k(L^d)=x_0(\alpha L^{d-k}x_1(L)\cdots
x_k(L)+\omega)=\]
\[\alpha (d-k)L^{d-(k+1)}x_0(L)x_1(L)\cdots x_k(L)+\]
\[\sum_j \alpha L^{d-k}x_1(L)\cdots x_0(x_j(L))\cdots x_k(L) +\]
\[\sum_i\alpha_i(d-(k-1))L^{d-k}x_0(L)x^i_1(L)\cdots x^i_{k-1}(L)+\]
\[\sum_i\sum_l\alpha_iL^{d-(k-1)}x^i_1(L)\cdots x_0(x^i_l(L))\cdots
x^i_{k-1}(L).\]
Let $z_j(L)=x_0(x_j(L))$ and $z^i_l(L)=x_0(x^i_l(L))$. Such elements
exist since $\lg/\lpl \otimes L\cong V$ as $P$-module.
Let 
\[ \omega= \sum_j\alpha L^{d-k}x_1(L)\cdots z_j(L)\cdots x_k(L)+\]
\[\sum_i \alpha_i(d-(k-1))L^{d-k}x_0(L)\cdots x^i_1(L)\cdots
x^i_{k-1}(L) +\]
\[\sum_i \sum_l \alpha_i L^{d-(k-1)}x^i_1(L)\cdots z^i_l(L)\cdots
x^i_{k-1}(L).\]
it follows $\omega \in L^{d-k}\otimes \sym^k(\lg/\lpl \otimes L)$.
Moreover
\[ x_0x_1\cdots x_k(L^d)=\tilde{\alpha} L^{d-(k+1)}x_0(L)\cdots
x_k(L)+\omega \]
where $\tilde{\alpha}=(d-k)\alpha$.
The Proposition is proved.
\end{proof}

\begin{theorem} \label{main} There is for all $1\leq l < d$ an isomorphism
\[ \U_l(\lg)L^d\cong L^{d-l}\otimes \sym^l(\lg/\lpl\otimes L) \]
of $P$-modules.
\end{theorem}
\begin{proof}  There are embeddings of $P$-modules
\[ \U_l(\lg)L^d \subseteq \sym^d(V) \]
and
\[ L^{d-l}\otimes \sym^l(\lg/\lpl \otimes L) \subseteq \sym^d(V).\]
Recall from \cite{flag} it follows
$dim_K(\U_l(\lg)L^d)=\binom{l+n}{n}$ where $dim_K(V)=n+1$.
Assume $z=x_1\cdots x_l(L^d)\in \U_l(\lg)L^d$. It follows from
Proposition \ref{envelop}
\[z=\alpha L^{d-l}x_1(L)\cdots x_l(L)+\omega\]
where 
\[ \omega \in L^{d-(l-1)}\sym^{l-1}(\lg/\lpl \otimes L)\subseteq
L^{d-l}\otimes \sym^l(\lg/\lpl \otimes L).\]
Since 
\[ \alpha L^{d-l}x_1(L)\cdots x_l(L)\in L^{d-l}\otimes \sym^l(\lg/\lpl
\otimes L) \]
it follows  $z\in L^{d-l}\otimes \sym^l(\lg/\lpl \otimes L) $ Hence we
get an inclusion of $P$-modules $\U_l(\lg)L^d \subseteq L^{d-l}\otimes
\sym^l(\lg/\lpl
\otimes L)$. Since
\[ dim_K(\U_l(\lg)L^d)=dim_K( L^{d-l}\otimes \sym^l(\lg/\lpl \otimes
L)  ) \]
the Theorem follows.
\end{proof}

\begin{corollary} \label{recov} There is for all $1\leq l <d$ an isomorphism
\[ \Pr^l_\p(\O_\p(d))(x) \cong (L^*)^{d-l}\otimes \sym^l(V^*) \]
of $P$-modules.
\end{corollary}
\begin{proof} There is by \cite{flag}, Theorem 3.10 an isomorphism
\[ \Pr^l_\p(\O_\p(d))(x)^*\cong \U_l(\lg)L^d \]
of $P$-modules. From this isomorphism and Theorem \ref{main} the
Corollary follows since 
\[ (L^{d-l}\otimes \sym^l(\lg/\lpl \otimes L))^*\cong
(L^*)^{d-l}\otimes \sym^l(V^*) \]
as $P$-modules.
\end{proof}

Note: Corollary \ref{recov} is proved in \cite{maa1} Theorem 2.4 using
more elementary techniques.

Let $Y=\Spec(K)$ and $\pi:\p(V^*)\rightarrow Y$ be the structure
morphism. Let $\p=\p(V^*)$.
Since $\sym^l(V^*)$ is a finite dimensional $\SL(V)$-module
it follows it is a free $\O_Y$-module with an
$\SL(V)$-linearization. It follows $\pi^*\sym^l(V^*)$ is a locally
free $\O_\p$-module with an $\SL(V)$-linearization since $\pi^*$
preserves the $\SL(V)$-linearization.

\begin{proposition} \label{slv} There is for all $1\leq l <d$ an isomorphism
\[ \Pr^l_\p(\O_\p(d))\cong \O_\p(d-l)\otimes \pi^*\sym^l(V^*) \]
of locally free $\O_\p$-modules with an $\SL(V)$-linearization.
\end{proposition}
\begin{proof} Let $P\subseteq \SL(V)$ be the subgroup fixing the line $L\in V$
There is an exact equivalence of categories
\begin{align}
&\label{equiv} \underline{mod}(P)\cong \ogpmod .
\end{align}
The $P$-module corresponding to $\O_\p(d-l)\otimes \pi^*\sym^l(V^*)$ is 
$(L)^{d-l}\otimes \sym^l(V^*)$. By the equivalence \ref{equiv} and
Corollary \ref{recov} we get an isomorphism
\[ \Pr^l_\p(\O_\p(d))\cong \O_\p(d-l)\otimes \pi^*\sym^l(V^*) \]
of locally free sheaves with $\SL(V)$-linearization and the
Proposition is proved.
\end{proof}

We get a formula for the splitting type of $\Pr^l_\p(\O_\p(d))$ on
projective space:

\begin{corollary} \label{split} There is for all $1\leq l < d$ an isomorphism
\[ \Pr^l_\p(\O(d))\cong \oplus^{\binom{n+l}{n}}\O_\p(d-l) \]
of locally free sheaves.
\end{corollary}
\begin{proof} The $P$-modules $\sym^l(V^*)$ corresponds to the free
  $\O_\p$-module $\oplus^{\binom{n+l}{n} } \O_\p$. 
The Corollary now follows from Proposition \ref{slv}.
\end{proof}

Let $X=\gr(m,m+n)$ and consider the $P$-modules
\[ L^{d-1}\otimes \sym^1(\lg/\lpl \otimes L) \subseteq \sym^d(\wedge^m
V) \]
and
\[ \U_1(\lg)L^d \subseteq \sym^d(\wedge^m V) .\]

\begin{proposition} \label{prop1} There is an isomorphism
\[ \U_1(\lg)L^d \cong L^{d-1}\otimes \sym^1(\lg/\lpl \otimes L) \]
of $P$-modules.
\end{proposition}
\begin{proof} Pick an element $x(L^d)=dL^{d-1}x(L)\in
  \U_1(\lg)L^d$. It follows $dL^{d-1}x(L)\in L^{d-1}\otimes
  \sym^1(\lg/\lpl \otimes L)$ hence there is an inclusion
\[ \U_1(\lg)L^d \subseteq L^{d-1}\otimes \sym^1(\lg/\lpl \otimes L) .\]
Let $L^{d-1}x(L)\in L^{d-1}\otimes \sym^1(\lg/\lpl \otimes L)$. It
follows
\[ L^{d-1}x(L)=\frac{1}{d}x(L^d)\in \U_1(\lg)L^d \]
hence there is an inclusion $L^{d-1}\otimes \sym^1(\lg/\lpl \otimes
L)$ and the Proposition is proved.
\end{proof}

\begin{corollary} \label{maincorr} There is an isomorphism
\[ \Pr^1_X(\O_X(d))(x)^*\cong L^{d-1}\otimes \sym^1(\lg/\lpl \otimes
L) \]
of $P$-modules.
\end{corollary}
\begin{proof} There is by \cite{flag}, Theorem 3.10 an isomorphism
\[ \Pr^1_X(\O_X(d))(x)^*\cong \U_1(\lg)L^d \]
of $P$-modules. The Corollary follows from this fact and
Proposition \ref{prop1}.
\end{proof}

Note: By \cite{maa10}, Example 5.12 there is a double complex
\[ \O_X(j)\otimes \H^i(X,\wedge^j \Pr^l_X(\O_X(d))^*) \]
of sheaves on $\p(W^*)$ where $W=\H^0(X,\O_X(d))$ and
$X=\gr(m,m+n)$. This double complex might give rise to a resolution of
the ideal sheaf of the $l$'th discriminant $D^l(\O_X(d))\subseteq
\p(W^*)$ of the line bundle $\O_X(d)$. By \cite{maa10}, Theorem 5.2 it
follows knowledge on the $P$-module structure of $\Pr^l_X(\O_X(d))$
gives information on the $\SL(V)$-module structure of the higher cohomology groups
$\H^i(X,\wedge^j \Pr^l_X(\O_X(d))^*)$ for all $i\geq 0$. This again gives information on the
dimension $h^i(X,\wedge^j \Pr^l_X(\O_X(d))^*)$. We get a description
of the locally free sheaf 
\[ \O_X(j)\otimes \H^i(X,\wedge^j \Pr^l_X(\O_X(d))^*) .\]
for all $i,j$.

\begin{example} \label{g24} Canonical filtration for the grassmannian $\gr(2,4)$.\end{example}
Consider the example where $m=n=2$ and $X=\gr(2,4)$. We get two
inclusions
\[ L^{d-2}\otimes \sym^2(\lg/\lpl \otimes L) \subseteq \sym^d(\wedge^2
V) \]
and
\[ \U_2(\lg)L^d \subseteq \sym^d(\wedge^2 V).\]
We may choose a basis for $\lp \subseteq \lg$ on the following form:
\[ \lp=\lpl \oplus L_x \]
where $L_x$ is the line spanned by the following vector $x$:
\[
x=\begin{pmatrix} 0 & 0 & 0 & 0 \\ 
                  0 & 1 & 0 & 0 \\
                  0 & 0 & -1 & 0 \\
                  0 & 0 & 0 & 0 
.\end{pmatrix}
\]
Let $\ln \subseteq \lg$ be the sub Lie algebra spanned by the
following vectors:
\[
x_1=
\begin{pmatrix} 0 & 0 & 0 & 0 \\
                  0 & 0 & 0 & 0 \\
                  1 & 0 & 0 & 0 \\
                  0 & 0 & 0 & 0
\end{pmatrix}
\]
\[
x_2=
\begin{pmatrix} 0 & 0 & 0 & 0 \\
                  0 & 0 & 0 & 0 \\
                  0 & 1 & 0 & 0 \\
                  0 & 0 & 0 & 0
\end{pmatrix}
\]
\[
x_3=
\begin{pmatrix} 0 & 0 & 0 & 0 \\
                  0 & 0 & 0 & 0 \\
                  0 & 0 & 0 & 0 \\
                  1 & 0 & 0 & 0
\end{pmatrix}
\]
and
\[
x_4=
\begin{pmatrix} 0 & 0 & 0 & 0 \\
                  0 & 0 & 0 & 0 \\
                  0 & 0 & 0 & 0 \\
                  0 & 1 & 0 & 0
\end{pmatrix}
.\]
Let $\tilde{\ln}$ be the vector space spanned by the vectors
$x_1,x_2,x_4,x_4$ and $x$.
It follows $\U_2(\lg)L^d=\U_2(\tilde{\ln})L^d\subseteq \sym^d(\wedge^2
V)$. The vector space $V$ has a basus $e_1,e_2,e_3$ and $e_4$. The
vector space $W$ has basis $e_1,e_2$. It follows $\wedge^2
W$ has a basis given by $e_1\wedge e_2=e[12]$ and $\wedge^2 V$ has basis
given by $e[12],e[13],e[14],e[23],e[24],e[34]$. By definition $L=e[12]$.
We get the following calculation:
\[ x_1(L)=-e[23], x_2(L)=e[13], x_3(L)=-e[24] \]
\[ x_4(L)=e[14], x(L)=e[12].\]
A basis for the $P$-module $L^{d-2}\otimes \sym^2(\lg/\lpl \otimes L)
$ are the following vectors:

\[ L^{d-2}x(L)x(L)=L^{d-2}e[12]^2 \]
\[L^{d-2}x_2(L)x(L)=L^{d-2}e[12]e[13]\]
\[ L^{d-2}x_4(L)x(L)=L^{d-2}e[12]e[14] \]
\[L^{d-2}x_1(L)x(L)=-L^{d-2}e[12]e[23] \]
\[ L^{d-2}x_3(L)x(L)=-L^{d-2}e[12]e[24] \]

\[ L^{d-2}x_2(L)x_2(L)=L^{d-2}e[13]^2  \]
\[ L^{d-2}x_2(L)x_4(L)=L^{d-2}e[13]e[14] \]
\[ L^{d-2}x_1(L)x_2(L)=-L^{d-2}e[13]e[23] \]
\[ L^{d-2}x_2(L)x_3(L)=-L^{d-2}e[13]e[24] \]

\[ L^{d-2}x_4(L)x_4(L)=-L^{d-2}e[14]^2  \]
\[ L^{d-2}x_1(L)x_4(L)=-L^{d-2}e[14]e[23] \]
\[ L^{d-2}x_3(L)x_4(L)=-L^{d-2}e[14]e[24] \]

\[ L^{d-2}x_1(L)x_1(L)=L^{d-2}e[23]^2e \]
\[ L^{d-2}x_1(L)x_3(L)=L^{d-2}e[23]e[24] \]

\[ L^{d-2}x_3(L)x_3(L)=L^{d-2}e[24]^2 \]

Let $a=d(d-1)$.
A basis for the $P$-module $\U_2(\lg)L^d=\U_2(\tilde{\ln})L^d$  are
the following vectors:
\[ x^2(L^d)=L^{d-2}e[12]^2 \]
\[ x_2x(L^d)=aL^{d-2}e[12]e[13]+dL^{d-1}e[13] \]
\[ x_4x(L^d)=aL^{d-2}e[12]e[14]+dL^{d-1}e[14]  \]
\[ x_1x(L^d)=aL^{d-2}e[12]e[23]-dL^{d-1}e[23] \]
\[ x_3x(L^d)= aL^{d-2}e[12]e[24]-dL^{d-1}e[24] \]

\[ x_2^2(L^d)=  aL^{d-2}e[13]^2 \]
\[ x_2x_4(L^d)= aL^{d-2}e[13]e[14] \]
\[ x_1x_2(L^d)= aL^{d-2}e[13]e[23] \]
\[ x_2x_3(l^D)= -aL^{d-2}e[13]e[24]-dL^{d-1}e[34] \]

\[ x_4^2(L^d)= L^{d-2}e[14]^2 \]
\[ x_1x_4(L^d)= -aL^{d-2}e[14]e[23]+dL^{d-1}e[34] \]
\[x_3x_4(L^d)= -aL^{d-2}e[14]e[24] \]

\[ x_1^2(L^d)=aL^{d-2}e[23]^2 \]
\[ x_1x_3(L^d)=aL^{d-2}e[23]e[24] \]

\[x_3^2(L^d)= aL^{d-2}e[24]^2 .\]

In the case where $W\subseteq V$ have dimensions $m$ and $m+n$ we get
embeddings of $P$-modules
\[ \U_l(\lg)L^d \subseteq \sym^d(\wedge^m V) \]
and 
\[ L^{d-l}\otimes \sym^l(\lg/\lpl \otimes L)\subseteq \sym^d(\wedge^m
V).\]
There is no equality
\[ \U_l(\lg)L^d=L^{d-l}\otimes \sym^l(\lg/\lpl \otimes L) \]
of $P$-modules as submodules of $\sym^d(\wedge^m V)$ in general as
Example \ref{g24} shows.

Since $\U_l(\lg)L^d$ and $L^{d-l}\otimes \sym^l(\lg/\lpl \otimes L)$
by Theorem \ref{main} and Proposition \ref{mainprop} 
are isomorphic when $m=1$ and $1\leq l<d$, have the same dimension
over $K$ and both have natural filtrations of $P$-modules we may
conjecture they are isomorphic as $P$-modules for all $m,n\geq 1$. 
Note: There is a canonical line $L^d\in \U_l(\lg)L^d$ for all
$l$. There is similarly a canonical line
\[ L^d\cong L^{d-l}\otimes \sym^l(\lp/\lpl\otimes L)\in L^{d-l}\otimes
\sym^l(\lg/\lpl \otimes L) .\]
Hence the two $P$-modules $\U_l(\lg)L^d$ and $L^{d-l}\otimes
\sym^l(\lg/\lpl \otimes L)$
look similar. 

In
general the $\SL(V)$-module $\sym^d(\wedge^m V)$ decompose
\[ \sym^d(\wedge^m V)\cong \oplus_i V_{\lambda_i}^{a_i} \]
where $V_{\lambda_i}$ are irreducible $\SL(V)$-modules and $a_i\geq 1$
are integers (see Proposition \ref{aut} for the situation of
$\gr(2,4)$). One may ask if there is a non-trivial automorphism
\[ \phi \in \aut_{\SL(V)}(\sym^d (\wedge^m V)) \]
with the property that the morphism
\[ \phi: \sym^d (\wedge^m V)\rightarrow \sym^d (\wedge^m V) \]
induce an isomorphism
\[ \tilde{\phi}:L^{d-l}\otimes \sym^l(\lg/\lpl \otimes L)
\rightarrow \U_l(\lg)L^d \]
of $P$-modules. In general the $\SL(V)$-module $\sym^d(\wedge^m V)$
has lots of automorphisms. When $m=2$ and $dim_K(V)=4$ it follows 
by Corollary \ref{aut} there is for every $d \geq 1$ an equality
\[ \aut_{\SL(V)}(\sym^d(\wedge^2 V) )= \prod_{i=0}^l K^*  \]
where $l=k$ if $d=2k$ or $d=2k+1$. For $m=n=2$ the $\SL(V)$-module
$\sym^d(\wedge^2 V)$ is by Proposition \ref{aut} multiplicity free.
The module $\sym^d(\wedge^m K^{m+n})$ is not multiplicity free in
general when $m,n>2$.

\section{Jet bundles and incidence complexes on the projective line}

In this section we construct a resolution by locally free sheaves of the ideal sheaf of
the $l$'th incidence scheme $I^l(\O_\p(d))\subseteq \p(W^*)\times \p$.
Here $\O_\p(d)$ is an invertible sheaf on the projective line $\p=\p^1$ and $W=\H^0(\p,\O_p(d))$.
There is on $Y=\p(W^*)\times \p^1$ a morphism $\phi(\O(d))$ of locally
free sheaves
\[ \phi(\O(d)):\O_{\p(W^*)}(-1)_Y\rightarrow \Pr^l(\O(d))_Y \]
Its zero scheme $Z(\phi(\O(d)))=I^l(\O(d))\subseteq Y$ is the $l$'th
incidence scheme of $\O(d)$. The Koszul complex of the morphism $\phi(\O(d))$
\[ 0\rightarrow \wedge^l\O(-1)_Y\otimes \Pr^l(\O(d))_Y^*\rightarrow
\cdots \rightarrow \wedge^2\O(-1)_Y\otimes \Pr^l(\O(d))_Y^*
\rightarrow \]
\[ \O(-1)_Y\otimes \Pr^l(\O(d))^*_Y \rightarrow \O_Y \rightarrow
\O_{I^l(\O(d))}\rightarrow 0 \]
- called the \emph{incidence complex of $\O(d)$} - is a
  resolution of the ideal sheaf of $I^l(\O(d))$. This follows from the
  fact that the ideal sheaf of $I^l(\O(d))$ is locally generated by a
  regular sequence.
We also calculate the higher direct images of the terms
\[ \O(-j)_Y\otimes \wedge^j \Pr^l(\O(d))^*_Y \]
appearing in the incidence complex.

The aim of the construction is to use it to construct a resolution of
the ideal sheaf of the discriminant $D^l(\O(d))$ where $\O(d)$ is a
line bundle on projective space or a grassmannian.

\begin{example} The Koszul complex of a map of modules.\end{example}

Let $A$ be an arbitrary commutative ring with unit and let
$\phi:E\rightarrow F$ be a map
$A$-modules.
Define the following map:
\[ d^1:E\otimes_A F^* \rightarrow A \]
by
\[ d^1(x\otimes f)=f(\phi(x)).\]
Let $I_\phi \subseteq A$ be the image of $d^1$. We let $I_\phi$ be
the \emph{ideal of the map $\phi$}.
Define the following map
\[ d^p:\wedge^pE\otimes F^* \rightarrow \wedge^{p-1}E\otimes F^* \]
by
\[ d^p(x_1\otimes f_1\wedge \cdots \wedge x_p\otimes
f_p)=\sum_{r=1}^p(-1)^{r-1}f_r(\phi(x_r))x_1\otimes f_1\wedge \cdots
\wedge \widetilde{x_r\otimes f_r}\wedge \cdots \wedge x_p\otimes
f_p.\]

\begin{lemma} \label{kcomplex} The following holds for all $p\geq 2$: $d^p\circ d^{p-1}=0$.
\end{lemma}
\begin{proof} We get
\[d^{p-1}d^p(x_1\otimes f_1\wedge \cdots \wedge x_p\otimes f_p)=\]
\[\sum_{r=1}^p(-1)^{r-1}f_r(\phi(x_r)) \]
\[\sum_{l\neq
  r}(-1)^{l-1}f_l(\phi(x_l))x_1\otimes f_1\wedge \cdots
\wedge \widetilde{x_l\otimes f_l} \wedge \cdots 
\wedge \widetilde{x_r\otimes f_r}\wedge \cdots \wedge x_p\otimes
f_p=0 \]
and the claim of the Lemma follows.
\end{proof}

Assume $E,F$ are locally free of finite rank and let $r=rk(E\otimes F^*)$.
We get a complex of locally free $A$-modules
\[ 0\rightarrow \wedge^rE\otimes F^*\rightarrow \cdots \rightarrow
\wedge^2E\otimes F^* \rightarrow E\otimes F^* \rightarrow A
\rightarrow A/I_\phi \rightarrow 0\]
called the \emph{Koszul complex of the map $\phi$}

\begin{example} \label{reg} The Koszul complex of a regular sequence.\end{example}
Let $\underline{x}=\{x_1,..,x_n\}$ be a regular sequence of elements
in $A$ and let $E=Ae$ be the free $A$-module on the element $e$. Let
$F=A\{e_1,..,e_n\}$ be a free rank $n$ module on $e_1,..,e_n$. Let $y_i=e_i^*$.
Define
\[ \phi:E\rightarrow F \]
by
\[ \phi(e)=x_1e_1+\cdots +x_ne_n.\]
Let $e\otimes y_i=z_i$. It follows
\[ d^p:\wedge^pE\otimes F^* \rightarrow \wedge^{p-1}E\otimes F^* \]
looks as follows:
\[ d^p(z_{i_1}\wedge \cdots \wedge z_{i_p})= \]
\[ \sum_{r=1}^p (-1)^{p-1}y_{i_r}(\phi(e))z_{i_1}\wedge
\cdots \wedge \widetilde{z_{i_r}} \wedge \cdots \wedge
z_{i_p}=\]
\[\sum_{r=1}^p(-1)^{r-1}x_{i_r}z_{i_1}\wedge
\cdots \wedge \widetilde{z_{i_r}} \wedge \cdots \wedge
z_{i_p}.\]
Hence the complex $\wedge^\bullet E\otimes F^*$ equals the Koszul
complex $K_\bullet(\underline{x})$ of the regular sequence
$\underline{x}$. It is an exact complex since $\underline{x}$ is a
regular sequence.

\begin{example} The Koszul complex of a morphism of locally free
  sheaves.\end{example}

The construction of the differential in the Koszul complex of a map of modules
is intrinsic, hence we may generalize to morphisms of 
locally free sheaves. Let $Y$ be an arbitrary  scheme and let
$\phi:\E\rightarrow \F$ be a map of locally free $\O_Y$-modules. 
Let 
\[ d^1:\E\otimes \F^* \rightarrow \O_Y \]
be defined locally by
\[ d^1(s\otimes v)=v(\phi(s)) .\]
Let $\I_\phi=Im(d^1)\subseteq \O_Y$ be the ideal sheaf defined by
$d^1$.  Since $\I_\phi$ is quasi coherent sheaf of ideals it follows the ideal sheaf
$\I_\phi$ corresponds to a subscheme $Z(\phi)\subseteq Y$ - the
\emph{zero scheme of $\phi$}.
Let
$U\subseteq Y$ be an open subset and define the following map:
\[ d^p:\wedge^p(\E\otimes \F^*)(U)\rightarrow \wedge^{p-1}(\E\otimes
\F^*)(U) \]
by
\[ d^p(s_1\otimes v_1\wedge \cdots \wedge s_p\otimes
v_p)=\sum_{r=1}^p(-1)^{r-1}v_r(\phi(s_r))s_1\otimes v_1\wedge \cdots
\wedge \widetilde{s_r\otimes v_r}\wedge \cdots \wedge s_p\otimes
v_p.\]
This gives a well defined map of locally free sheaves since we have
not chosen a basis for the module $\wedge^p(\E\otimes \F^*)(U)$ to
give a definition. 
By Lemma \ref{kcomplex} it follows $d^p\circ d^{p+1}=0$ for all $p\geq 1$
hence we get a complex of locally free sheaves.
The sequence of maps of locally free sheaves
\[ 0\rightarrow \wedge^r\E\otimes \F^*\rightarrow \cdots \rightarrow
\wedge^2\E\otimes \F^* \rightarrow \E\otimes \F^* \rightarrow  \O_Y\rightarrow \O_{Z(\phi)}\rightarrow 0 \]
is called the \emph{Koszul complex of the map $\phi$}. Here $r=rk(\E\otimes \F^*)$.

\begin{example} Koszul complexes and local complete intersections. \end{example}

Assume $\phi:\L\rightarrow \F$ is a map of locally free $\O_Y$-modules
where $\L$ is a line bundle. Let $Z(\phi)\subseteq Y$ be the subscheme
defined by $\phi$ - the zero scheme of $\phi$.
Let $r=rk(\F)$. Choose an open affine cover
$U_i$ of $Y$ where $\F$ and $\L$ trivialize, i.e
\[ \F(U_i)=\O(U_i)\{f_{i1},..,f_{ir} \} \]
and 
\[ \L(U_i)=\O(U_i)e_i .\]
Let $\O(U_i)=A_i$, $L_i=\L(U_i)$ and $F_i=\F(U_i)$.
Assume the image
\[ \phi(U_i):L_i\rightarrow F_i \]
has
\[ \phi(U_i)(e_i)=x_{i1}f_{i1}+\cdots +x_{ir}f_{ir} \]
where $\{x_{i1},..,x_{ir}\}\subseteq A_i$ is a regular sequence. Let
$I_i=\underline{x_i}=\{x_{i1},..,x_{ir}\}$. It follows from Example
\ref{reg} the Koszul complex
\[ 0\rightarrow \wedge^r(L_i\otimes F_i^*)\rightarrow \cdots
\rightarrow \wedge^2(L_i\otimes F_i^*) \rightarrow L_i\otimes F_i^* \]
\[ \rightarrow A_i\rightarrow A_i/I_i \rightarrow
0\]
is a resolution of the ideal $I_i$ since $I_i$ is generated by a
regular sequence. The complex $\wedge^\bullet L_i\otimes F_i^*$ is
isomorphic to the Koszul complex $K_\bullet(\underline{x_i})$ on the
regular sequence $\underline{x_i}$.
It follows the global complex
\[ 0\rightarrow \L^{\otimes r}\wedge^r  \F^*\rightarrow \cdots \rightarrow
\L^{\otimes 2}\wedge^2\F^* \rightarrow \L\otimes \F^* \rightarrow
\O_Y\rightarrow \O_{Z(\phi)}\rightarrow 0 \]
is a resolution of the ideal sheaf $\I_{Z(\phi)}$ of $Z(\phi)\subseteq
Y$ since it is locally isomorphic to the Koszul complex
$K_\bullet(\underline{x_i})$ for all $i$. 

Since the ideal $I_i$ is
generated by a regular sequence of lenght $r$ it follows
$dim(A_i/I_i)=dim(A_i)-r$. If $Y$ is irreducible of dimension $d$ it
follows $Z(\phi)\subseteq Y$ is a local complete intersection of
dimension $d-r$.

\begin{example} The incidence complex of $\O(d)$ on the projective line. \end{example}

Let $\p=\p^1_K$ where $K$ is a field of characteristic zero and let
$\O(d)\in \Pic(\p)=\Z$ be a line bundle where $d\in \Z$. Let 
\[ W=\H^0(\p,\O(d))=K\{e_0,..,e_d \} \]
where $e_i=x_0^{d-i}x_1^i$. Let $y_i=e_i^*$. Let $Y=\p(W^*)\times \p$
and consider the following diagram
\[
\diagram   Y \rto^p \dto^q & \p \dto^\pi \\
            \p(W^*) \rto^\pi & \Spec(K) 
\enddiagram.
\]
There is a sequence of locally free $\O_Y$-modules
\[ \O_{\p(W^*)}(-1)_Y\rightarrow \H^0(\p.\O(d))\otimes \O_Y
\rightarrow^{T^l_Y} \Pr^l(\O(d))_Y \]
and let $\phi(\O(d))$ be the composed map
\begin{align}
&\label{composed} \phi(\O(d)):\O_{\p(W^*)}(-1)_Y\rightarrow \Pr^l(\O(d))_Y .
\end{align}
It follows by \cite{maa10} the zero scheme $Z(\phi(\O(d)))$ equals the
incidence scheme $I^l(\O(d))$ of the line bundle $\O(d)$.
By definition $\p(W^*)=\Proj(K[y_0,..,y_d])$ where $y_i=e_i^*$. It has
an open cover on the following form: $D(y_i)=\Spec(K[u_0,..,u_d])$
where we let $u_j=\frac{y_j}{y_i}$. Let $y_j/y_j=1$. Let
\[ F(t)=u_0+u_1t+\cdots +u_dt^d\in K[u_0,..,u_d,t] .\]
Restrict the map \ref{composed} to the open set $U_{i0}=D(y_i)\times
D(x_0)\subseteq Y$. 
We get the following two maps of modules:
\[ \alpha: \O_{\p(W^*)}(-1)|_{U_{i0}}\rightarrow \O_{U_{i0}}\otimes
\H^0(\p,\O(d)) \]
\[\alpha: K[y_i,t]\frac{1}{y_i}\rightarrow K[u_i,t]\otimes
K\{e_0,..,e_d\} \]
defined by
\[ \alpha(1/y_i)=\sum_{k=0}^d u_k\otimes e_k =\sum_{k=0}^d u_k\otimes
x_0^{d-k}x_1^k=
\sum_{k=0}^d u_k\otimes t^kx_0^d.\]
We get the map
\[ T^l_{U_{i0}}:\O_{U_{i0}}\otimes \H^0(\p,\O(d)) \rightarrow
\Pr^l(\O(d))|_{U_{i0}} \]
defined by
\[ T^l(1\otimes x_0^{d-i}x_1^i)=T^l(1\otimes t^ix_0^d)=(t+dt)^i\otimes
x_0^d.\]
The composed map
\[ \phi(\O(d))_{U_{i0}}:K[u_i,t]\frac{1}{y_i}\rightarrow
K[u_i,t]\{dt^j\otimes x_0^d\} \]
is the map
\[ \phi(\O(d))(\frac{1}{y_i})=\sum_{k=0}^du_k(t+dt)^k\otimes
x_0^d=\]
\[ \sum_{k=0}^l \frac{F^{(k)}(t)}{k!}dt^k \otimes x_0^d\in
K[u_i,t]\{1\otimes x_0^d,..,dt^l\otimes x_0^d\} .\]

Let $U_{i1}=D(y_i)\times D(x_1)\subseteq Y$ and let
$\frac{x_0}{x_1}=s$. Let 
\[ G(s)=u_d+u_{d-1}s+u_{d-2}s^2+\cdots +u_0s^d\in K[u_0,..,u_d,s].\]
Restrict the map \ref{composed} to the open set $U_{i1}$

We get the following two maps of modules:
\[ \alpha: \O_{\p(W^*)}(-1)|_{U_{i1}}\rightarrow \O_{U_{i1}}\otimes
\H^0(\p,\O(d)) \]
\[\alpha: K[y_i,s]\frac{1}{y_i}\rightarrow K[u_i,s]\otimes
K\{e_0,..,e_d\} \]
defined by
\[ \alpha(1/y_i)=\sum_{k=0}^d u_k\otimes e_k =\sum_{k=0}^d u_k\otimes
x_0^{d-k}x_1^k=
\sum_{k=0}^d u_k\otimes s^{d-k}x_1^d.\]
We get the map
\[ T^l_{U_{i1}}:\O_{U_{i1}}\otimes \H^0(\p,\O(d)) \rightarrow
\Pr^l(\O(d))|_{U_{i1}} \]
defined by
\[ T^l(1\otimes x_0^{d-i}x_1^i)=T^l(1\otimes s^{d-i}x_1^d)=(s+ds)^{d-i}\otimes
x_1^d.\]
The composed map
\[ \phi(\O(d))_{U_{i1}}:K[u_i,s]\frac{1}{y_i}\rightarrow
K[u_i,s]\{ds^j\otimes x_1^d\} \]
is the map
\[ \phi(\O(d))(\frac{1}{y_i})=\sum_{k=0}^du_{d-k}(s+ds)^k\otimes
x_1^d=\]
\[ \sum_{k=0}^l \frac{G^{(k)}(s)}{k!}ds^k \otimes x_1^d\in
K[u_i,s]\{1\otimes x_1^d,..,ds^l\otimes x_1^d\} .\]
It follows the ideal sheaf $\I_{I^l(\O(d))}$ of $I^l(\O(d))$ is
generated by
\[ \{ \frac{F^{(l)}(t)}{l!}, \frac{F^{(l-1)}(t)}{(l-1)!},.., F(t)  \}
\]
on $U_{i0}$ and by
\[ \{ \frac{G^{(l)}(s)}{l!}, \frac{G^{(l-1)}(s)}{(l-1)!},..,G(s) \}   \]
on $U_{i1}$. Let $z_i=\frac{F^{(i)}(t)}{(i)!}$ and $w_i=\frac{G^{(i)}(s)}{(i)!}$
for $i=0,..,l$.
\begin{lemma} \label{lemma} Assume $B$ is a commutative ring of characteristic zero
  and let
\[ f(t)=a_0+a_1t+\cdots +a_dt^d \in B[t] \]
be an arbitrary degree $d$ polynomial with $a_d\neq 0$. Let
$f^{(i)}(t)$ denote the formal derivative with respect to $t$.
It follows
\[ \frac{f^{(k)}(t)}{l!}=\sum_{i=k}^d \binom{i}{k}a_it^{i-k}.\]
\end{lemma}
\begin{proof} The proof is by induction. 
\end{proof}

\begin{lemma} The sequence $\{z_l,..,z_0\}$ is a regular sequence in
  $K[u_i,t]$. The sequence $\{w_l,..,w_0\}$ is a regular sequence in $K[u_i,s]$.
\end{lemma}
\begin{proof}  Let $z_i=\frac{F^{(i)}(t)}{i!}$ and
  $w_j=\frac{G^{(j)}}{j!}$.
Assume $l<i$ and consider the sequence $z_l,z_{l-1},..,z_0\subseteq
A[t]=K[u_0,..,u_d][t]$.
Since $A[t]$ is a domain it follows $z_l$ is a non zero divisor in
$A[t]$. We see from Lemma \ref{lemma}
\[ A[t]/w_l\cong K[u_0,..,u_{l-1},u_{l+1},..,u_d,t] \]
which is a domain, hence $w_{l-1}$ is a non zero divisor in
$A[t]/w_l$. By induction it follows $z_l,..,z_0$ is a regular sequence
in $A[t]$. Assume $i\leq l$.
It follows the sequence $z_l,..,z_{i+1}$ is a regular sequence in
$A[t]$. We see  from Lemma \ref{lemma}  $z_i$ is non zero in 
\[ A[t]/(z_l,..,z_{i+1})=K[u_0,..,u_i,u_{l+1},..,u_d,t] \]
and $K[u_0,..,u_i,u_{l+1},..,u_d,t]$ is a domain. It follows $z_i$ is
a non zero divisor in $A[t]/(z_l,..,z_{i+1})$. It follows $z_l,..,z_0$
is a regular sequence in $A[t]$ and the claim follows.
A similar argument proves $w_l,..,w_0$ is a regular sequence in $A[s]$
and the Lemma is proved.
\end{proof}

One may prove using similar methods for any permutation $\sigma \in S_{l+1}$ the sequences
\[ z_{\sigma(l)},..,z_{\sigma(0)} \]
and
\[ w_{\sigma(l)},..,w_{\sigma(0)} \]
are regular sequences.

It follows the ideal sheaf $\I_{I^l(\O(d))}$ is locally generated by a
regular sequence.

The morphism
\[ \phi(\O(d)):\O_{\p(W^*)}(-1)_Y\rightarrow \Pr^l(\O(d))_Y \]
gives by Example \ref{reg} rise to a Koszul complex 
\[ \wedge^\bullet \O_{\p(W^*)}(-1)\otimes \Pr^l(\O(d))^*_Y \]
of locally free sheaves of $Y=\p(W^*)\times \p^1$.

\begin{definition}
Let the complex
\begin{align}
&\label{incidence} 0\rightarrow \wedge^l\O(-1)_Y\otimes \Pr^l(\O(d))_Y^*\rightarrow
\cdots \rightarrow \wedge^2\O(-1)_Y\otimes \Pr^l(\O(d))_Y^*
\rightarrow 
\end{align}
\[ \O(-1)_Y\otimes \Pr^l(\O(d))^*_Y \rightarrow \O_Y \rightarrow
\O_{I^l(\O(d))}\rightarrow 0 \]
be the  \emph{incidence complex of $\O(d)$}.
\end{definition}

Since the ideal sheaf of $I^l(\O(d))$ by the discussion above is locally generated by a regular
sequence it follows from Example \ref{reg} the complex \ref{incidence} is a resolution.

In \cite{maa1}, Theorem 5.10 one calculates the higer direct images
\[ \R^iq_*(\wedge^j\O(-1)_Y\otimes \Pr^l(\O(d))^*_Y) \]
for all $i,j$. We get the following calculations:

Let $V=K\{e_0,e_1\}$ and $\p=\p(V^*)$. Let
$W=\H^0(\p,\O(d))=\sym^d(V^*)$ and consider the diagram
\[ 
\diagram Y=\p(W^*)\times \p \rto^p \dto^q & \p \dto^\pi \\
           \p(W^*) \rto^\pi & \Spec(K) 
\enddiagram.
\]
By the results of this paper it follows there is an isomorphism
\[ \Pr^l_\p(\O(d))\cong \O_\p(d-l)\otimes \pi^* \sym^l(V^*) \]
a sheaves with an $\SL(V)$-linearization. 
We get
\[ \wedge^j\Pr^l_\p(\O_\p(d))\cong \O_\p(j(d-l))\otimes \pi^*\wedge^j
\sym^l(V^*) .\]
By the equivariant projection formula for higher direct images we get
\[ \R^i q_*(\wedge^j\O(-1)_Y\otimes
\Pr^l(\O(d))^*_Y)\cong \O_{\p(W^*)}(-j)\otimes
\H^i(\p,\wedge^j\Pr^l_\p(\O_\p(d))^*). \]
Let  
\[ \pi:\p\rightarrow \Spec(K) .\]
It follows
\[ \wedge^j\Pr^l_\p(\O_\p(d))^*\cong \O_\p(j(l-d))\otimes
\pi^*\wedge^j \sym^l(V^*) .\]
We get
\[ \H^i(\p,\wedge^j \Pr^l(\O(d))^*)\cong
\R^i\pi_*(\pi^*(\wedge^j\sym^l(V))\otimes \O_\p(j(l-d)))\cong \]
\[ \wedge^j(\sym^l(V))\otimes H^i(\p, \O_\p(j(l-d))).\]
We get the following Theorem:

\begin{theorem} The following holds:
\[\R^ip_*(\O(-j)\otimes \wedge^j\Pr^l(\O(d))^*)=0 \text{ if $i=0$ or
  $i=1$ and $j(d-l)<2$.}\]
\[ \R^1p_*(\O(-j)\otimes \wedge^j\Pr^l(\O(d))^*)=\O(-j)\otimes
\sym^{j(d-l)-2}(V)\otimes \wedge^j\sym^l(V) \]
if $j(d-l)\geq 2$.
\end{theorem}
\begin{proof} The proof follows from the calculation of the
  equivariant cohomology of line bundles on projective space (see \cite{jantzen}).
\end{proof}

Hence we have complete control on the sheaf
\[ \R^i q_*(\wedge^j\O(-1)_Y\otimes \Pr^l(\O(d))^*_Y) \]
on the projective line and projective space for all $i,j$. Using the techniques
introduced in this paper one may describe resolutions of incidence
schemes $I^l(\O(d))$ on more general grassmannians and flag varieties.
The hope is we may be able to construct resolutions of the ideal sheaf of $D^l(\O(d))$ using
indicence resolutions in a more general situation.

Note: In \cite{lascoux} resolutions of ideal sheaves of
determinantal schemes are studied and much is known on such
resolutions. In \cite{maa10} it is proved $D^1(\O(d))$ is a
determinantal scheme for any $d\geq 2$ on the projective line
$\p^1$. Assume $\L\in \Pic^G(G/P)$ is a $G$-linearized linebundle, $G$ a semi
simple linear algebraic group and $P$ a parabolic subgroup. If one can
prove $D^l(\L)$ is a determinantal scheme we get two approaches to the
study of resolutions of ideal sheaves of discriminants: One using jet
bundles and incidence schemes, another one using determinantal schemes.

\section{Appendix A: Automorphisms of representations}

Let $W \subseteq V$ be vectorspaces of dimension two and four over the
field $K$.
Consider the subgroup $P\subseteq G=\SL(V)$ where $P$ is the parabolic
subgroup of elements fixing $W$. It follows $\pi: G\rightarrow G/P=\gr
(2,4)$ is a principal $P$-bundle. Let $\lg=Lie(G)$ and $\lp=Lie(P)$ be the Lie
algebras of $G$ and $P$. In this section we study the decomposition
into irreducibles and automorphisms of some $G$-modules. We also study
some $\Pss$-modules where $\Pss$ is  the semi-simplification of $P$.
It follows $\Pss$ equals $\SL(2)\times \SL(2)$. Since $\lp
\subseteq \lg$ is a
$P$-sub module it follows the quotient $\lg/\lp$ is a $P$-module
hence a $\Pss$ module. We may apply the theory of highest weights
since $\Pss=\SL(2)\times \SL(2)$ is a semi simple algebraic group.

\begin{proposition} \label{prop1} The following hold: There is an isomorphism of
  $\SL(2)\times \SL(2)$-modules
\begin{align}
\sym^k(\lg / \lp)=& \oplus_{i=0}^n \sym^{2i+m}(W^*)\otimes
\sym^{2i+m}(V/W).
\end{align}
for all $k\geq 1$.
Here $(n,m)=(\frac{k}{2},0)$ if $k=2n$ and $(n,m)=(\frac{k-1}{2},1)$ if $k=2n+1$.
\end{proposition}
\begin{proof}  Recall the canonical isomorphism from Lemma \ref{canonical}
\[ \lg/\lp\cong \Hom(W,V/W)\cong W^*\otimes V/W \]
of $P$-modules. It follows 
\[ \sym^k(\lg/\lp)\cong \sym^k(W^*\otimes V/W) \]
and its decomposition into irreducible $\SL(2)\times \SL(2)$-modules
can be done using well known formulas (see \cite{fulton}). Alternatively one may compute
its highest weight vectors and highest weights explicitly using the construction from
Section 5.
\end{proof}

Let $i:G/P\rightarrow \p(\wedge^2 V^*)=\p$ be the Plucker embedding and
let $\O_{G/P}(1)=i^*\O_{\p}(1)$ be tautological
line bundle on $G/P$ and let $\O_{G/P}(d)=\O_{G/P}(1)^{\otimes d}$.
It follows from the Borel-Weil-Bott Theorem
$\H^0(\gr, \O_{\gr}(d))$ is an irreducible $\SL(V)$-module.
Let $V$ have basis $e_1,e_2,e_3,e_4$ and let $\wedge^2 V$ have basis
$e_{ij}$ for $1\leq i < j \leq 4$, with $e_{ij}=e_i\wedge e_j$.
Consider the element $f\in \sym^2(\wedge ^2 V)$ where
\[ f=e_{12}e_{34}-e_{13}e_{24}+e_{14}e_{23} .\]
One checks $f$ is a highest weight vector for $\SL(V)$ with
highest weight $0$, hence it defines the unique trivial character of
$\SL(V)$.
Its dual 
\[ f^*=x_{12}x_{34}-x_{13}x_{24}+x_{14}x_{23}\in \sym^2(\wedge^2 V^*) \] 
is the defining equation for $\gr=G/P$ as closed subscheme of $\p(\wedge^2 V^*)$.

\begin{proposition} \label{aut} The following hold: there is an
  isomorphism of
$\SL(V)$-modules
\begin{align}
\sym^d(\wedge^2 V)=& \oplus_{i=0}^l \H^0(\gr, \O_{\gr}(d-2i) )^*,
\end{align}
where $l=k$ if $d=2k$ or $d=2k+1$.
\end{proposition}
\begin{proof} The result is proved using the theory of highest
  weights. There is a split exact sequence of $\SL(V)$-modules
\[ 0 \rightarrow f^*\sym^{d-2}(\wedge^2 V^*) \rightarrow
\sym^d(\wedge^2 V^*) \rightarrow \H^0(\gr, \O_{\gr}(d) ) \rightarrow 0
.\]
Dualize this sequence to get the split exact sequence
\[ 0 \rightarrow f\sym^{d-2}(\wedge^2 V) \rightarrow
\sym^d(\wedge^2 V) \rightarrow Q_d \rightarrow 0 .\]
where $Q_d=\H^0(\gr, \O_{\gr}(d) )^*$. Since $f$ is the trivial
character it
follows there is an isomorphism
\[ f\sym^d(\wedge^2 V) \cong \sym^d(\wedge ^2 V) \]
of $\SL(V)$-modules. By the Borel-Weil-Bott Theorem it
follows  $Q_d$ is an irreducible $\SL(V)$-module.
If $d=2k$ we get by induction the equality
\[ \sym^d(\wedge^2 V^*)=Q_{d}\oplus Q_{d-2}\oplus \cdots \oplus Q_2
\oplus Q_0,\]
and the claim of the Proposition is proved in the case where
$d=2k$.
The claim when $d=2k+1$ follows by a similar argument and the
Proposition is proved.
\end{proof}

\begin{corollary} Let $\E=\oplus_{i=0}^l \O_\gr(2i-d)$ where
  $l=k$ if $d=2k$ or $d=2k+1$.
It follows
\[ \H^0(\gr, \E)\cong \sym^d(\wedge^2 V^*) \]
as $\SL(V)$-module.
\end{corollary}
\begin{proof} We get by Proposition \ref{aut} isomorphisms of $\SL(V)$-modules
\[ \H^0(\gr, \E)\cong\H^0(\gr, \oplus_{i=0}^l \O_\gr(d-2i))\cong\]
\[\oplus_{i=0}^l
\H^0(\gr,\O_\gr(d-2i))\cong\sym^d(\wedge^2V)^*\cong\sym^d(\wedge^2
V^*) \]
and the Corollary is proved.
\end{proof}

\begin{corollary} \label{aut} There is for every $d \geq 1$ an equality
\[ \aut_{\SL(V)}(\sym^d(\wedge^2 V) )= \prod_{i=0}^l K^*  \]
where $l=k$ if $d=2k$ or $d=2k+1$.
\end{corollary}
\begin{proof} This follows from Proposition \ref{aut} and the
  Borel-Weil-Bott theorem (BWB). From the BWB theorem it follows 
  $\H^0(\gr, \O_{\gr}(d))^*$ is an
  irreducible $\SL(V)$-module for all $d\geq 1$. From this and
  Proposition \ref{aut} the claim of the Corollary follows.
\end{proof}

Hence the $\SL(V)$-module $\sym^d(\wedge^2V)$ is a multiplicity free
$\SL(V)$-module for all $d\geq 1$. This is not true in general for
$\sym^d(\wedge^m K^{m+n})$ when $m,n> 2$.

In general if $\mathbb{S}_\lambda$ and $\mathbb{S}_\mu$ are two
\emph{Schur-Weyl modules} (see \cite{fulton}) there is a decomposition
\[ \mathbb{S}_\lambda (\mathbb{S}_\mu(V))\cong \oplus_i  V_{\lambda_i}
\]
where $V_{\lambda_i}$ is an irreducible $\SL(V)$-module for all $i$.
It is an open problem to calculate this decomposition for
two arbitrary partitions $\lambda $ and $\mu$.

\section{Appendix B: The Cauchy formula}

We include in this section an elementary discussion of the Cauchy
formula using multilinear algebra. Let
$W\subseteq V$ be vector spaces of dimension $m$ and $m+n$ over $K$ and let
$P\subseteq \SL(V)$ be the subgroup fixing $W$. Let $\lg=Lie(G)$ and
$\lp=Lie(P)$. There is a canonical isomorphism
\[ \lg/\lp \cong \Hom(W, V/W) \]
of $P$-modules, hence the elements of $\lg/\lp$ may be interpreted
as
linear maps. The symmetric power $\sym^k(\lg/\lp)=\sym^k(\Hom(W,V/W))$
is a $P$-module hence a $P_{semi}=\SL(m)\times \SL(n)$-module and we
want to give an
explicit construction of its highest weight vectors as $\Pss$-module.

\begin{proposition} \label{determinant} Let $U=K^m$. There is a
  canonical map of
  $\SL(V)$-modules
\[ \wedge^m(U^*)\otimes \wedge^mU \rightarrow \sym^m(\Hom(U,U)) \]
defined by
\[ x_1\wedge \cdots \wedge x_m\otimes e_1\wedge \cdots \wedge e_m
\rightarrow
\begin{vmatrix} x_1\otimes e_1 & x_1\otimes e_2 & \cdots & x_1\otimes
  e_m \\
x_2\otimes e_1 & x_2\otimes e_2 & \cdots & x_2\otimes
  e_m \\
x_m\otimes e_1 & x_m\otimes e_2 & \cdots & x_m\otimes
  e_m .
\end{vmatrix}
\]
Here $e_1,..,e_m$ is a basis for $U$ and $x_1,..,x_m$ is a basis for
$U^*$.
\end{proposition}
\begin{proof}The proof is left to the reader as an exercise.
\end{proof}

Note: in Proposition \ref{determinant} the element $x_i\otimes e_j$ is
an element of $U^*\otimes U=\Hom(U,U)$. Hence the determinant
\[
\begin{vmatrix} x_1\otimes e_1 & x_1\otimes e_2 & \cdots & x_1\otimes
  e_m \\
x_2\otimes e_1 & x_2\otimes e_2 & \cdots & x_2\otimes
  e_m \\
x_m\otimes e_1 & x_m\otimes e_2 & \cdots & x_m\otimes
  e_m .
\end{vmatrix}
\]
may be interpreted as a polynomial of degree m in the elements
$x_i\otimes e_j$, hence it is an element of $\sym^m(\Hom(U,U))$.

Let $B\subseteq \SL(m,K)\times \SL(n,K)\subseteq \SL(V)=\SL(m+n,K)$
be the following subgroup:
$B$ consists of matrices with determinant one of the form
\[
\begin{pmatrix} U_1 & 0 \\
                0 & U_2
\end{pmatrix}
\]
where
\[U_1=
\begin{pmatrix} a_{11} & 0 & \cdots & 0 \\
                a_{21} & a_{22} & \cdots & 0 \\
 \vdots & \vdots  & \cdots & \vdots \\
 a_{m1} & a_{m2}  & \cdots & a_{mm} \\
\end{pmatrix}
\]
and
\[ U_2
\begin{pmatrix} b_{11} & 0 & \cdots & 0 \\
b_{21}  & b_{22} & \cdots & 0 \\
\vdots & \vdots & \cdots & \vdots \\
b_{n1}  & 0 & b_{n2} \cdots & b_{nn}
\end{pmatrix}.
\]

Let $T$ be a $B$-module and $v\in T$ a vector with the
property that for all $x\in B$ it follows 
\[ xv=\lambda(x)v \]
where $\lambda \in \Hom(B,K^*)$ is a character of $B$. It follows
$v$ is a highest weight vector for $T$ as $\SL(m,K)\times
\SL(n,K)$-module.
The group $B   \subseteq \SL(V)$ defines  filtrations of $W$ and $V/W$
as follows: Let $W$ have basis $e_1,..,e_m$ and $V$ have basis
$e_1,..,e_m,f_1,..,f_n$.
Let $W_1=\{e_m\}$, $W_2=\{e_m, e_{m-1}\}$ and
\[ W_i=\{e_m,..e_{m-i+1} \} .\]
It follows we get a filtration
\[ 0=W_0\subseteq W_1\subseteq \cdots \subseteq W_{m-1}=W \]
of the vector space $W$.
Let
\[ U_j=W_{m-1}\cup \{f_n,..,f_{n-j+1} \} \]
and let $V_i=(V/W)/U_{n-i}$. We get a surjection
\[ V/W \rightarrow V_i \]
for $i=1,..,n-1$.
It follows $dim W_i=dim V_i=d_i$ for all $i$.
Let $x:W\rightarrow V/W$ be a linear map of vector spaces.
We get an induced map
\[ x_i:W_i\rightarrow V_i \]
wich is a square $d_i$ matrix for all $i$.
Let $g\in B$ be the element

\[
\begin{pmatrix} G_1 & 0 \\
                0 & G_2
\end{pmatrix}
\]

where
\[G_1=
\begin{pmatrix} a_{1} & 0 & \cdots & 0 \\
                * & a_{2} & \cdots & 0 \\
 \vdots & \vdots  & \cdots & \vdots \\
 * & *   & \cdots & a_{m} \\
\end{pmatrix}
\]
and
\[G_2
\begin{pmatrix} b_{1} & 0 & \cdots & 0 \\
*  & b_{2} & \cdots & 0 \\
\vdots & \vdots & \cdots & \vdots \\
* & 0 &  \cdots & b_{n}
\end{pmatrix}.
\]

The $i'th$ wedge product
\[ |x_i|=\wedge^i x_i\in \Hom(\wedge^i W_i,\wedge^i
V_i)=\wedge^i(W_i^*)\otimes \wedge^iV_i \]
may be viewed as an element in
\[ |x_i|\in \sym^i(\Hom(W_i,V_i))\subseteq \sym^i(\Hom(W,V/W)) \]
via Proposition \ref{determinant}.

\begin{proposition} The following formula holds:
\[ g|x_i|=\frac{b_1\cdots b_i}{a_{m-i+1}\cdots
  a_m}|x_i|=\lambda(g)|x_i| \]
for all $g\in B$. Here $\lambda(g)=\frac{b_1\cdots
  b_i}{a_{m-i+1}\cdots
  a_m}$ is a character $\lambda \in \Hom(B,K^*)$.
\end{proposition}
\begin{proof} The proof is left to the reader as an exercise.
\end{proof}
Hence the $i$'th determinant $|x_i|\in \sym^i(\Hom(W,V/W))$  is a
highest weight vector for the $\SL(m)\times \SL(n)$-module
$\sym^i(\Hom(W,V/W))$.
By the results of \cite{brion} it follows  the vectors
$x_0^{d_0}x_1^{d_1}\cdots x_i^{d_i}$ with $\sum id_i=k$ are all
highest weight vectors for the module 
\[ \sym^k(\Hom(W,V/W))\cong \sym^k(W^*\otimes V/W).\]

\textbf{Acknowledgements}. The author thanks Michel Brion, Alexei
Roudakov  and an anonymous referee for comments on the contents of this paper.

\end{document}